\documentclass[12pt]{amsart}
\usepackage{amssymb}
\usepackage{amsfonts}
\usepackage{latexsym}

\input xy
\xyoption{matrix}\xyoption{arrow}\xyoption{curve}
\def\arrow{\ar}

\long\def\ignore#1{#1}

\def\uloopr#1{\arrow@'{@+{[0,0]+(-4,5)} @+{[0,0]+(0,10)} @+{[0,0]+(4,5)}}
  ^{#1}}
\def\dloopr#1{\arrow@'{@+{[0,0]+(-4,-5)} @+{[0,0]+(0,-10)} @+{[0,0]+(4,-5)}}
  _{#1}}
\def\rloopd#1{\arrow@'{@+{[0,0]+(5,4)} @+{[0,0]+(10,0)} @+{[0,0]+(5,-4)}}
  ^{#1}}
\def\lloopd#1{\arrow@'{@+{[0,0]+(-5,4)} @+{[0,0]+(-10,0)} @+{[0,0]+(-5,-4)}}
  _{#1}}

\def\gldim{\mbox{gl\,dim\,}}
\def\fd{\mbox{fin\,dim\,}}
\def\Fd{\mbox{Fin\,dim\,}}
\def\Kdim{\mbox{K\,dim\,}}
\def\pdim{\mbox{p\,dim\,}}
\def\lamod{\Lambda\mbox{-mod}}
\def\pinf{{\mathcal P}^\infty}
\def\calO{{\mathcal O}}
\def\A{{\mathcal A}}
\def\C{{\mathcal C}}
\def\len{\mbox{length\,}}
\def\laMod{\Lambda\mbox{-Mod}}
\def\Soc{\mbox{Soc\,}}
\def\rank{\mbox{rank}}
\def\card{\mbox{card\,}}
\def\Hom{\mbox{Hom}}
\def\la{{\Lambda}}

\begin{document}

\title[The finitistic dimension conjectures -- a tale of 3.5 decades]{The finitistic dimension conjectures -- a tale of 3.5 decades}
\author[Birge Zimmermann Huisgen]{Birge Zimmermann Huisgen}
\address{
Department of Mathematics \\
University of California\\
Santa Barbara, CA 93106 \\
U.S.A.
}
\email{birge@math.ucsb.edu}

\begin{abstract} We review the history of several homological conjectures,
both from a chronological and a methodological point of view.
\end{abstract}

\keywords{Finite dimensional algebra, finitistic dimension,
projective resolution, projective dimension, syzygy, contravariant finiteness.}

\maketitle

\section{Introduction}

Since the last five years have brought fresh initiatives in connection with
the Finitistic Dimension Conjectures, providing not only new results but, at
the same time, new homological methods, a synopsis of the recent developments
appears helpful in keeping stock.

After some brief remarks on the origins of homological algebra, we will
present the Finitistic Dimension Conjectures publicized by Bass in 1960, as
well as their connections with other problems concerning the homology and
structure theory of finite dimensional algebras and classical orders. We will
give an overview over results to date, followed by a non-technical outline of
the ideas involved. The methods applied in this context can be roughly
categorized as follows; here $\la$ will stand for a finite dimensional
algebra. 

(1) Detecting repetitions in the sequence $\Omega^1(M)$,
$\Omega^2(M)$, $\Omega^3(M)$, $\dots$ of syzygies of a $\la$-module $M$.

(2) Passing numerical data along projective resolutions.

(3) Finding the structure of second or higher syzygies.

(4) Determining when the category $\pinf$ of finitely generated $\la$-modules
of finite projective dimension is contravariantly finite in $\lamod$.

\section{A bit of history}

I should really go back farther than I intend to, since the idea of using
chain complexes to measure the distance of a given object from a situation
which is considered ideal goes back to topology where, for instance, singular
chain complexes were used to measure the discrepancy of a given space from
simple connectedness. Around the middle of this century, algebraists started to
adapt this idea, in the first place for use in commutative algebra, in
particular, for a better understanding of the rings that arise as coordinate
rings of algebraic varieties. 

Officially, homological algebra came to life in the early 1950's, with
Cartan, Eilenberg, Nakayama, Auslander, Buchsbaum, Serre, and Nagata being
among the principal instigators. However, there is one result going back to
the last century which is, by all means, a homological one. In 1890,
Hilbert's celebrated work ``\"Uber die Theorie der algebraischen Formen''
appeared \cite{Hil}, and in it a theorem which we now label
\medskip

\noindent{\bf  Hilbert's Syzygy Theorem (1890).} \emph{If $k$ is a field, then ${\rm
gl\,dim\,} k[X_1,\dots,X_n] =n$.}
\medskip

Starting his argument with the wry comment ``Der Beweis ist nicht
m\"uhelos'', Hilbert actually computed by brute force sets of relations for
the iterated syzygies of a finitely generated module over $k[X_1,\dots,X_n]$,
and arrived at the conclusion that the $n$-th generation of relations is
necessarily trivial. The ten pages of heavy-duty computations of Hilbert's
original proof leave us somewhat relieved that mathematicians like Emmy Noether
provided a conceptual framework for the ``modern'' reasoning. On the side: it
was also Hilbert who introduced the term ``syzygy'' into algebra in its
present meaning (the term had previously occurred in a paper of Sylvester in
1853). Going back to the Greek roots ``syn--zygon'', meaning ``yoked
together'', the term was first used by astronomers referring to
configurations of three celestial bodies in a straight line.

The second milestone I want to mention is the
\medskip

\noindent{\bf  Auslander-Buchsbaum-Serre Theorem (1955/56).} \emph{If $\,V$ is an algebraic
variety over an algebraically closed field and $R$ its coordinate ring, then
the global dimension of $R$ is finite if and only if $V$ is smooth. Moreover,
in the smooth case, ${\rm gl\,dim\,} R =\dim V$.}
\medskip

\noindent{\bf  Corollary (Auslander-Buchsbaum, 1959).} \emph{If $\,V$ is smooth, then all
localizations of $R$ are unique factorization domains.}
\medskip

For details, see \cite{Ser}, \cite{AB56}, \cite{AB57}, \cite{AB59}. The
question of unique factorization had been pursued for a considerable period,
and part of the credit should go to Nagata \cite{Nag} and Zariski, who
independently reduced the problem to varieties of dimension 3. It is this
reduced situation which Auslander and Buchsbaum successfully tackled with
homological methods \cite{AB59}.

As alluded to earlier, the general idea behind introducing homological
dimensions, also into the representation theory of noncommutative rings, was to
find a measure for the deviation of a given module category from the
``ideal'' categories arising in the Artin-Wedderburn situation which is
characterized by the projectivity of all objects. This
enterprise of measuring the complexity of a module category in terms of the
projective dimensions of its objects and that of a ring
by means of its global dimension was clearly very successful in that, in
many cases, these dimensions do provide a very effective measure. Let me,
however, remind you of an example where the global dimension does not at all
live up to expectations.
\medskip

{\bf Example.} Let $R= k[X]/(X^2)$ where $k$ is a field. Then

$\bullet$ The projective dimension of the $R$-module $R/(X)$ is infinite; in
fact, all syzygies of $R/(X)$ are isomorphic to $R/(X)$. Hence, $\gldim R$ is
infinite. On the other hand, the category of
$R$-modules is extremely simplistic, namely,

$\bullet$ Each $R$-module is of the form $\bigl[ \bigoplus$ copies of
$R\bigr]$
$\oplus$ $\bigl[\bigoplus$ copies of $R/(X)\bigr]$.
\medskip

In this and many other examples, a far more accurate measure of the
complexity of the module category is provided by the finitistic dimensions of
which I will remind you next:
\begin{eqnarray*}
\mbox{(l)}\fd R &= &\sup \{ \pdim M \mid M \mbox{\ a fin. gen. left\ }
R\mbox{-module with\ } \pdim M <\infty \}\\
\mbox{(l)}\Fd R &= &\sup \{ \pdim M \mid M \mbox{\ an arbitrary left\ }
R\mbox{-module with\ } \pdim M <\infty \} 
\end{eqnarray*}
Start by observing that in the previous example, $\Fd R =0$, which is more in
line with the simplicity of the corresponding module category.

The little finitistic dimension, $\fd$, occurred implicitly already in the
proof of the Auslander-Buchsbaum-Serre Theorem, while Kaplansky suggested
also studying the analogous invariant obtained by waiving the restriction to
finitely generated modules. In exploring these two new invariants, the first
two natural questions appear to be as follows: Do they coincide? Are they
always finite? Both questions were promptly answered in the negative for
noetherian rings, even in the commutative noetherian situation.
\medskip

$\bullet$ If $R$ is commutative noetherian local, then $\fd R =\mbox{depth\,}
R$. The argument is implicit in Auslander's and Buchsbaum's paper \cite{AB57}.
 So in particular, a commutative noetherian local ring has coinciding big and
little finitistic dimensions if and only if it is Cohen-Macaulay. 

$\bullet$ If $R$ is commutative noetherian, then $\Fd R =\Kdim R$. The
inequality $\ge$ was established by Bass in 1962 \cite{B62}; the other was
filled in by Gruson and Raynaud in 1971 \cite{GR71}. In particular, $\Fd R$
may be infinite: indeed, commutative noetherian rings of
infinite Krull dimension are well known to exist, Nagata having supplied the
first examples. 
\medskip

However, no examples settling these questions were available in the
non\-com\-mu\-ta\-tive artinian situation. (Of course, commutative artinian
rings are uninteresting in this context, since all the invariants in question
are zero in that case.) This led to the following conjectures, which were
publicized as ``problems'' by Bass in 1960 and later -- restricted to finite
dimensional algebras -- promoted to the status of conjectures.
\medskip

{\bf Finitistic Dimension Conjectures.} Let $\Lambda$ be a finite
dimensional algebra over a field $k$.

(I) $\fd\Lambda =\Fd\Lambda$;

(II) $\fd\Lambda<\infty$.
\medskip

I will start with a very compressed overview of results pertaining directly to
the conjectures. Next, I will indicate how strongly the conjectures are
interwoven with other representation theoretic problems, and finally, I will
revisit existing results, again sketching them in fairly broad strokes, but
this time from a methodological standpoint.

\section{Synopsis of results}

Throughout, $\Lambda$ will be a finite dimensional algebra over a field $k$,
with Jacobson radical $J$.
\medskip

$\bullet$ If $\fd\Lambda =0$, then Conjecture I holds, i.e., $\Fd\Lambda =0$
as well. This easy observation is embedded in a far stronger result of Bass
\cite{B60}. To verify it, suppose that $\Fd\Lambda >0$, which guarantees the
existence of a left $\Lambda$-module $M$ with $\pdim M=1$. Let $f :
P\rightarrow M$ be a projective cover. Then the kernel of $f$ is nonzero
and projective, and contained in the radical of $P$. In particular, this
forces (a copy of) a principal left ideal $\Lambda e$, where $e$ is a primitive
idempotent, to be contained in $JP$. Clearly, $\Lambda e$ is then contained
in $JP_0$, where $P_0$ is a finitely generated projective submodule of $P$, and
$P_0/\Lambda e$ is a finitely generated module of projective dimension 1.

$\bullet$ If $J^2=0$, the second Conjecture is valid; more precisely, we then
have $\Fd\Lambda \le 1+ \sup\{ \pdim S\mid S\in\lamod \mbox{\ simple of
finite projective dimension} \}$. This is due to
Mochizuki \cite{M65}. For a proof, it suffices to observe that the first
syzygy of any $\Lambda$-module is semisimple.

$\bullet$ A strengthening of the previous result: If the projective dimension
of $J^2$, viewed as a {\it right} $\Lambda$-module, is finite, Conjecture II
is true. Small derived this from Mochizuki's result, with the following
change of rings trick \cite{S69}, applied to $I=J^2$: If $R$ is any ring, $I$
a nilpotent two-sided ideal, and $M$ a left $R$-module with $\mbox{Tor}_m^R
(R/I,M)=0$ for all $m>0$, then $\mbox{p\,dim}_R\, M= \mbox{p\,dim}_{R/I}\,
M/IM$.
\medskip

The next series of events in this connection occurred only after a long time
gap.
\medskip

$\bullet$ Suppose that $\Lambda$ is a monomial relation algebra, i.e.,
$\Lambda= k\Gamma/I$, where $\Gamma$ is a quiver and $I$ an ideal generated
by certain paths in $k\Gamma$. Then Conjecture II holds.

The first to prove this were Green, Kirkman, and Kuzmanovich in 1991
\cite{GKK}. Simplified proofs and explicit bounds on $\Fd\Lambda$ were given
by Igusa-Zacharia \cite{IZ}, Cibils \cite{Ci}, and the author \cite{ZHpre}. We
will comment on the ideas involved later.

$\bullet$ Conjecture II holds whenever $J^3=0$. It is even sufficient to
replace the hypothesis by the weaker condition $\pdim J^3_\Lambda <\infty$.

This was proved by Green and the author in 1991 \cite{G-ZH}. Simplified
arguments are due to Fuller-Saorin \cite{FuSa} and Igusa-Todorov \cite{IT}.
Moreover, the following generalization appeared in two instalments, the first
of which is due to Dr\"axler-Happel \cite{DH}, and the second to Wang
\cite{Wan}: If
$J^{2l+1}=0$ and
$\Lambda/J^l$ has finite representation type, then $\fd\Lambda<\infty$.

$\bullet$ Denote by $\pinf$ the full subcategory of $\lamod$ the objects
of which are the modules of finite projective dimension. If $\pinf$ is
contravariantly finite in $\lamod$ (the definition will follow in Section 8
below), Conjecture II is true. This was shown by Auslander and Reiten, again
in 1991
\cite{AR91}.

$\bullet$ Conjecture I fails, even for monomial relation algebras. In fact,
for each $n\ge 2$, there exists a monomial relation algebra $\Lambda$ such
that $\fd\Lambda =n$ and $\Fd\Lambda =n+1$. Examples were given by the
author in 1992 \cite{ZHinvent}. A simplified proof for the underlying theory
was developed by Butler \cite{But}.
\medskip

In order to widen the horizon, I would like to point out a few of the numerous
connections between the Finitistic Dimension Conjectures and other problems.
Moreover, I would like to include some immediately contingent questions in the
discussion.

\section{Implications and related problems}

$\bullet$ If $\fd\Lambda<\infty$, the Nakayama Conjecture holds for
$\Lambda$. This latter conjecture asserts that a finite dimensional algebra
$\Lambda$ of infinite dominant dimension (i.e., possessing a minimal
injective resolution $0\rightarrow {}_\Lambda\Lambda\rightarrow Q_0\rightarrow
Q_1\rightarrow \cdots$ in which all the terms $Q_i$ are projective) is
quasi-Frobenius. This implication was pointed out by Tachikawa in
\cite{Tac}. I include the easy argument: Denote the homomorphisms in the above
resolution by $f_i : Q_i\rightarrow Q_{i+1}$, and suppose that $\fd\Lambda
=d<\infty$. Assuming that all of the $Q_i$ are projective, the kernel of
$f_i$ has projective resolution

$$0\rightarrow \Lambda\rightarrow Q_0\rightarrow \cdots \rightarrow
Q_{i-1}\rightarrow \ker f_i\rightarrow 0.$$

Combining this with the fact that $\pdim \ker f_{d+1} \le d$, we see that
$\Lambda$ is a direct summand of
$Q_0$ and thus injective.

$\bullet$ In case $\fd\Lambda<\infty$, finiteness of the left injective
dimension of $\Lambda$ implies finiteness of the right injective dimension.
(It is well known that these dimensions coincide when both are finite.) This
was pointed out by Auslander.

$\bullet$ Problem: Find explicit bounds on the finitistic dimensions of
$\Lambda$ in terms of $\dim_k \Lambda$.

Recall the following non-constructive result due to Schofield \cite{Sch}: Given
a field
$k$, there exists a map $f : \bf N\rightarrow \bf N$ such that, for each
finite dimensional $k$-algebra $\Lambda$, either $\gldim \Lambda=\infty$ or
else $\gldim\Lambda \le f(\dim_k \Lambda)$. It was subsequently observed by
Jensen and Lenzing \cite{JLbook} that an analogous statement holds if the
global dimension of $\Lambda$ is replaced by $\sup \{ \pdim M\mid M\in
\lamod,\ \len M\le b\}$, where $b$ is a fixed positive integer. Clearly this
gives rise to the hope for a corresponding result for finitistic dimensions, at
least for those known to be finite. Closely related is the problem of finding
test classes of modules on which the finitistic dimension is known to be
attained. In contrast to the global dimension, the finitistic dimensions are in
general not attained on the cyclic modules, not even on the $b$-generated
modules for any prefixed integer $b$ (see Section 7 below). 

To underline the usefulness
of such bounds, we briefly describe an old problem on classical orders which
can be translated into a question concerning a rather specialized, but still
poorly understood, class of finite dimensional algebras.

$\bullet$ Problem: Establish bounds on the global dimension of classical orders
over discrete valuation rings. Let $\calO\subseteq M_n(K)$
be a classical order over
$D$, where
$D$ is a discrete valuation ring with quotient field $K$; moreover, denote by
$\pi$ a uniformizing parameter of $D$. In 1970, Tarsy conjectured that the
finiteness of the global dimension of $\calO$ implies $\gldim\calO\le n-1$, where
$n$ is the matrix size of the order \cite{Tar}. This conjecture may have been
based merely on a paucity of examples; nonetheless, it appears to have had a
stimulating effect on the area. In \cite{VaJ2}, V. A. Jategaonkar observed
that, up to isomorphism, the number of orders in $M_n(K)$ of finite global
dimension is finite for any $n$, whence a bounding function of $n$ on the
finite global dimensions arising inside $M_n(K)$ does exist. Moreover, positive
results relating to the conjecture were obtained for specialized classes of
orders by her \cite{VaJ1,VaJ2}, Kirkman-Kuzmanovich
\cite{KK}, and Fujita \cite{Fuj2}. In particular, the conjecture
received some credence in the situation where $\calO$ is tiled, meaning that
$\calO$ contains a full set of $n$ orthogonal idempotents. However, in
\cite{Fuj2}, Fujita exhibited a class of tiled classical orders in $M_n(K)$,
$n\ge 6$, having global dimension $n$, thus refuting the conjecture. Recently,
Jansen and Odenthal constructed examples showing that a replacement of
$n-1$ by $n$ will not save the conjecture; in fact, they obtained tiled
classical orders
$\calO_n\subseteq M_n(K)$, $n\ge8$ and even, such that $\gldim\calO_n =2n-8$
\cite{JO}. So this leaves the problem of suggesting another plausible bound
in terms of the matrix size. A suggestion for tackling this question was
given by Green-Kirkman-Kuzmanovich in \cite{GKK}, where they show the
following homological connection between a classical order $\calO$ and the
finite dimensional algebra $\Lambda =\calO/\pi\calO$ over the field $D/(\pi)$;
namely,
$\fd\calO =1+\fd\Lambda$. Moreover, Kirkman and Kuzmanovich observed that
$\gldim\Lambda =\infty$ whenever
$n\ge2$ and
$\calO$ is tiled. The reason why this approach hs a chance of being more viable
lies in the fact that deleting
$\pi\calO$ gives rise to a considerable simplification of the situation; in
other words, the homology of $\Lambda$ should be more transparent than that
of $\calO$.

$\bullet$ Question: What is the structure of modules in $\lamod$ having finite
projective dimension, in contrast to the structure of those of infinite
projective dimension? One of the major problems encountered in connection
with the Finitistic Dimension Conjectures lies in the fact that, short of
computing projective resolutions, one has no means of recognizing modules of
finite projective dimension. A promising key to such structure theorems was
provided by Auslander-Reiten in \cite{AR91} for the case where $\pinf$ is
contravariantly finite in $\lamod$ (see Section 8 below).

$\bullet$ Problem: Develop algorithms for computing or estimating homological
dimensions. In \cite{AnGr}, Anick-Green laid the foundations for an algorithm
to construct projective resolutions; this work was later combined with an
adapted theory of Gr\"obner bases by Farkas-Feustel-Green \cite{FFG} and
implemented on the computer by Feustel and Green. A low-effort algorithm for
obtaining the finitistic dimensions of monomial relation algebras, up to an
error of 1, was presented by the author in
\cite{ZHpre}. Furthermore, a computer program for the computation of global
dimensions of classical orders was developed at the University of Stuttgart
under the direction of K. W. Roggenkamp.  

\medskip
At this point, we will supplement the sketch of the existing results given
above, and discuss the principal ideas behind the arguments.

\section{Method 1: Repetitions in the structure of syzygies}

This method was used in proving the first generation of results for monomial
relation algebras. 

The idea is to take advantage of repetitions likely to occur in the sequence
of syzygies of a module. Recall the trivial example given to motivate
the introduction of the finitistic dimensions; in that example, all syzygies
recurred ad infinitum. While this example is particularly simplistic, the
phenomenon as such is not at all atypical. However, usually it does not occur
in quite as clean-cut a form.

To measure repetitiveness in the structure of syzygies, we introduce the
following {\it repetition index} of a module $M$ in $\laMod$:
\begin{eqnarray*}
\rho(M) &= &\inf \{i\in{\bf N} \mid \mbox{each non-projective
indecomposable summand of\ } \Omega^i(M)\\ 
 &&\hskip0.75truein \mbox{\ occurs as a summand of\ } \Omega^j(M)
\mbox{\ for infinitely many\ } j\}.
\end{eqnarray*}
In case the above set is empty, we set $\rho(M)=\infty$. This concept is not
altogether new; without introducing the repetition index, Jans \cite{Jan}
referred to modules $M$ with the property that $\mbox{add}(\bigoplus_{i\ge1}
\Omega^i(M))$ has finite representation type as modules having an ``ultimately
closed'' projective resolution; clearly such modules have finite repetition
index. Another slightly different concept measuring repetition is due to
Fuller and Wang \cite{FuWa}: they say that the projective resolution of a
module $M$ ``has a strongly redundant image from an integer $n$'' in case
each indecomposable direct summand of $\Omega^n(M)$ occurs as a summand of
some $\Omega^j(M)$ with $j>n$; this number $n$ may differ from the repetition
index by 1.
\medskip

$\bullet$ The following facts were implicitly proved by Igusa and Zacharia in
\cite{IZ}:

(a) If $\rho =\rho((\Lambda/J)_\Lambda)$, then $\rm{(l)}\Fd \Lambda \le \rho$.

(b) If $\Lambda$ is a monomial relation algebra, then
$\rho((\Lambda/J)_\Lambda) \le \dim_k J$.

Remarks: (i) The inequality under (a) remains true whenever $\rho
=\rho(E_\Lambda)$, where $\Soc E$ contains all simple right
$\Lambda$-modules \cite{GoHu}. Interestingly, injective cogenerators usually
yield better bounds on $\Fd \Lambda$ than $\Lambda/J$. On the side, we mention
that Fuller and Wang \cite{FuWa} extended a version of (a) to noetherian
rings.

(ii) If $\Lambda =\calO/\pi\calO$ for a tiled classical order $\calO$ as above, then
finiteness of the global dimension of $\calO$ implies that all $\Lambda$-modules
have finite repetition index; this was observed by Goodearl and the author
\cite{GoHu}. We conjecture that the repetition indices of injective
cogenerators for Mod-$\la$ are always finite when $\la =\calO/\pi\calO$. If
confirmed, this would show that the little finitistic dimensions of tiled
classical orders are necessarily finite.
\medskip

{\bf Example for Method 1}. Let $D= k[[\pi]]$, where $k$ is a field and $\pi$
an indeterminate, and let $K= k((\pi))$. If we define $\calO$ to be the
following subring of $M_5(K)$
$$\calO = \left( \begin{matrix} D&D&D&D&D \\
\pi&D&D&D&D \\ \pi^2&\pi&D&D&D \\ \pi^2&\pi^2&\pi&D&\pi \\
\pi^3&\pi^2&\pi&\pi&D \end{matrix} \right)$$
then $\Lambda =\calO/\pi\calO$ has the quiver
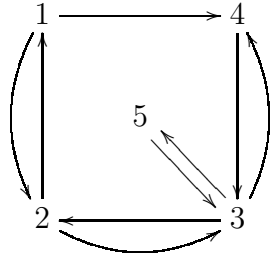
\ignore{
$$\xyoption{matrix}\xyoption{arrow}\xyoption{curve}
\xy\xymatrix{
1 \arrow[rr] \arrow@/_1pc/[dd] &&4 \arrow[dd]\\
 &5 \arrow[dr]<-0.5ex>\\
2 \arrow[uu] \arrow@/_1pc/[rr] &&3 \arrow[ll] \arrow@/_1pc/[uu]
\arrow[ul]<-0.5ex>
}\endxy$$}
\vskip0.1truein
\noindent and indecomposable left projectives $\Lambda e_1,\dots, \Lambda e_5$
with the following graphs:
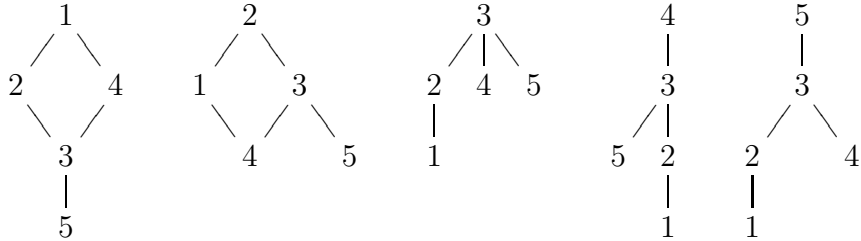
\ignore{
$$\xyoption{matrix}\xyoption{arrow}\xymatrixrowsep{1pc}\xymatrixcolsep{0.5pc}
\xy\xymatrix{
 &1 \arrow@{-}[dl] \arrow@{-}[dr] &&&&2 \arrow@{-}[dl] \arrow@{-}[dr] &&&&&3
\arrow@{-}[dl] \arrow@{-}[d] \arrow@{-}[dr] &&&&4 \arrow@{-}[d] &&&5
\arrow@{-}[d]\\
2 \arrow@{-}[dr] &&4 \arrow@{-}[dl] &&1 \arrow@{-}[dr] &&3 \arrow@{-}[dl]
\arrow@{-}[dr] &&&2 \arrow@{-}[d] &4 &5 &&&3 \arrow@{-}[dl] \arrow@{-}[d]
&&&3 \arrow@{-}[dl] \arrow@{-}[dr]\\
 &3 \arrow@{-}[d] &&&&4 &&5 &&1 &&&&5 &2 \arrow@{-}[d] &&2 \arrow@{-}[d] &&4\\
 &5 &&&&& &&&&& &&&1 &&1
}\endxy$$}
\noindent See \cite{ZHpre,BuHu} for an introduction to these graphs;
similar
graphical representations of modules were developed by Alperin \cite{Alp} and
Fuller \cite{Ful}. Here we mention only that the
entries in the
$l$-th row of vertices in the graph of a module $M$  give the indices of
the simple direct summands of $J^{l-1}M/J^lM$. It is readily checked that
the syzygies of
$S_1 =\Lambda e_1/Je_1$ have the following graphs:
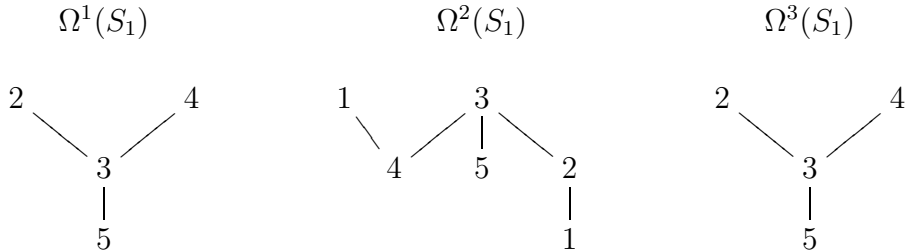
\ignore{
$$\xyoption{matrix}\xyoption{arrow}\xymatrixrowsep{1pc}\xymatrixcolsep{0.5pc}
\xy\xymatrix{
 &\Omega^1(S_1) &&&& &&&\Omega^2(S_1) &&& &&&\Omega^3(S_1)\\
2 \arrow@{-}[dr] &&4 \arrow@{-}[dl] &&&&1 \arrow@{-}[dr] &&3
\arrow@{-}[dl] \arrow@{-}[d] \arrow@{-}[dr] &&&&&2 \arrow@{-}[dr] &&4
\arrow@{-}[dl]\\
 &3 \arrow@{-}[d] &&& &&&4 &5 &2
\arrow@{-}[d] &&&&&3 \arrow@{-}[d]\\
 &5 &&&& &&&&1 &&&&&5
}\endxy$$}\noindent which shows that $\rho(S_1)=1$. We leave it as an easy exercise for
the reader to check that $\rho(S_2\oplus \dots\oplus S_5)=2$. This implies that
$$\fd \calO =1+\fd \Lambda \le 3.$$
(Here, in fact, the global dimension of $\calO$ is 3.)

\section{Method 2: Passing numerical data along projective resolutions}

This method applies to the case where $J^3=0$ as well as to monomial relation
algebras (see \cite{G-ZH,ZHbounds,ZHpre}). 

The underlying idea is roughly this: Let $M$ be a finitely generated
left $\Lambda$-module which is contained in the radical of a projective
module. (Note that all syzygies of finitely generated modules are of this
form.) The goal is to find a package of numerical data for $M$ which
determines whether $M$ is zero, and allows to compute the analogous package
of data for $\Omega^1(M)$. The numerical invariants I have in mind are
multiplicities of the simple modules in the various radical layers of $M$,
which can be conveniently arranged in a suitable group $\rm{Mat}\,\bf Z$ of
integer matrices. Suppose, for the moment, that for each module $M$ as above,
we have found a suitable matrix $[M]\in \rm{Mat}\,\bf Z$, together with a
$\bf Z$-linear map
$L :
\rm{Mat}\,\bf Z\rightarrow \rm{Mat}\,\bf Z$ such that
$$[M]=0 \Longleftrightarrow M=0$$
$$L([M]) =[\Omega^1(M)].$$
Then, clearly, the projective dimension of $M$ is less than or equal to $m$
if and only if $L^{m+1}([M]) =0$, and hence $\fd\Lambda \le \rank_{\bf Z}
(\rm{Mat}\,\bf Z)$.

The easiest instance where this method works to advantage is the case where
$J^3=0$ and all the simple left modules, $S_1,\dots,S_n$, have infinite
projective dimension. In that case, we choose $\rm{Mat}\,\bf Z$ to be the
group of all $2\times n$ matrices over $\bf Z$ and $[M]$ to be the matrix
listing the multiplicities of the $S_i$ in $M/JM$ and in $JM$ as entries of
the first and second row, respectively. Then we can easily find a linear map
$L$ as above, which tells us that
$\fd\Lambda
\le 2n$ in this case.
 Examples readily show that this bound is sharp. If some of the simples $S_i$
have finite projective dimension, a similar choice of matrices still
yields a bound on
$\fd\Lambda$, but not in terms of the vector space dimension of $\Lambda$; in
fact, the left finitistic dimension of $\Lambda$ is bounded above by
$$1+\sup\{ \pdim S_i\mid \pdim S_i<\infty\} +2\cdot \card\{ S_i\mid \pdim S_i
=\infty\}.$$
 (See
\cite{G-ZH}.)

Considering several matrix groups simultaneously, however, does lead to a
general bound in terms of the vector space dimension for algebras $\Lambda$
with vanishing radical cube \cite{ZHbounds}. Namely: 

\medskip
$\bullet$ If $J^3=0$ and $n$ is again the
number of isomorphism types of simple $\Lambda$-modules, then
$$\fd\Lambda \le n^2+1 \le (\dim_k \Lambda)^2+1.$$
\medskip

{\bf Example for Method 2}. Let $\Gamma$ be the following quiver:
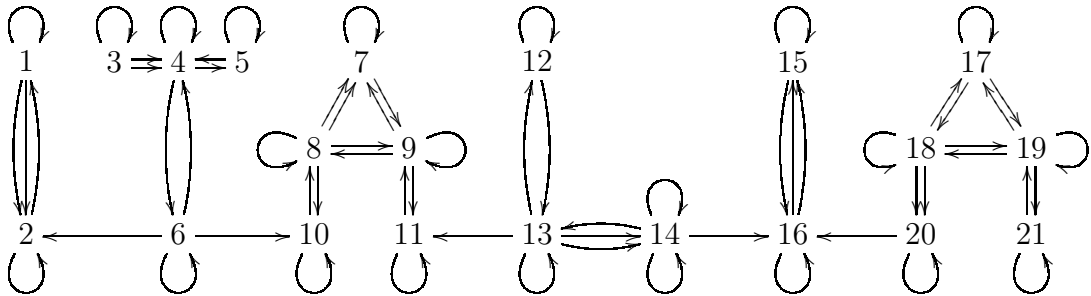
\ignore{
$$\xyoption{matrix}\xyoption{arrow}\xyoption{curve}
\xymatrixrowsep{1.5pc}\xymatrixcolsep{0.18pc}
\xy\xymatrix{
1 \uloopr{} \arrow[dd] \arrow@/_0.4pc/[dd] &&&3 \uloopr{} \arrow[rr]<0.3ex>
\arrow[rr]<-0.3ex> &&4 \uloopr{} \arrow[rr]<-0.3ex> \arrow@/_0.4pc/[dd] &&5
\uloopr{} \arrow[ll]<-0.3ex> &&&7 \uloopr{} \arrow[dr]<0.3ex> &&&&&12
\uloopr{}
\arrow@/^0.4pc/[dd] &&&& &&&&15 \uloopr{} \arrow@/_0.4pc/[dd] &&&&&17
\uloopr{}
\arrow[dl]<0.3ex> \arrow[dr]<0.3ex>\\
 &&&&& &&&&8 \lloopd{} \arrow[ur]<0.3ex> \arrow[ur]<-0.3ex>
\arrow[rr]<0.3ex>
\arrow[d]<0.3ex> &&9 \rloopd{} \arrow[ul]<0.3ex> \arrow[ll]<0.3ex>
\arrow[d]<0.3ex> &&&& &&&& &&&& &&&&18 \lloopd{} \arrow[ur]<0.3ex>
\arrow[rr]<0.3ex>
\arrow[d]<-0.3ex> \arrow[d]<0.3ex> &&19 \rloopd{} \arrow[ul]<0.3ex>
\arrow[ll]<0.3ex> \arrow[d]<0.3ex>\\
2 \dloopr{} \arrow@/_0.4pc/[uu] &&&&&6 \dloopr{}
\arrow[0,-5] \arrow@/_0.4pc/[uu] \arrow[0,4] &&&&10 \dloopr{}
\arrow[u]<0.3ex> &&11
\dloopr{} \arrow[u]<0.3ex> &&&&13 \dloopr{} \arrow[0,-4] \arrow[0,4]
\arrow@/_0.4pc/[0,4]
\arrow@/^0.4pc/[uu] &&&&14 \uloopr{} \dloopr{} \arrow@/_0.4pc/[0,-4]
\arrow[0,4] &&&&16
\dloopr{} \arrow[uu] \arrow@/_0.4pc/[uu] &&&&20 \dloopr{} \arrow[0,-4] &&21
\dloopr{} \arrow[u]<0.3ex> 
}\endxy$$}
\noindent Whenever $I\subseteq k\Gamma$ is an ideal containing all paths of
length 3,
$$\fd(k\Gamma/I) \le 42.$$

\section{Method 3: Finding the structure of second or higher syzygies}

This method provides the theory in the background of the above-mentioned
examples refuting the first Finitistic Dimension Conjecture.

``Usually'', the submodules of projectives which arise as $m$-th syzygies get
progressively simpler as $m$ grows. The key to results derived from such
simplifications lies in the observation that, if one completely
understands
$m$-th syzygies, one can determine the finitistic dimensions up to an error of
$m-1$. Monomial relation algebras provide an instance where this method
succeeds with
$m=2$. The following was proved in \cite{ZHpre}:\medskip

$\bullet$ If $\Lambda= k\Gamma/I$ is a monomial relation algebra and $M$ an
arbitrary left
$\Lambda$-module, then $\Omega^2(M)$ is isomorphic to a direct sum of copies
of principal left ideals $\Lambda p$ for suitable paths $p$ in $k\Gamma
\setminus I$ of positive length.
\medskip

The original proof of this result is somewhat rough going; those who prefer
it smoother we refer to \cite{But}. This
structure of second syzygies permits to compute both the big and the little
finitistic dimension of $\Lambda$ up to an error of 1. Namely, if
$$s= \sup \{ \pdim \Lambda q \mid q \mbox{\ a path of length\ } \ge 1 \mbox{\
with\ } \pdim \Lambda q <\infty\},$$
and $s= -1$ if this set is empty, then
$$\fd\Lambda,\ \Fd\Lambda \in [s+1,s+2].$$
There is a straightforward graphical method for computing $s$, which
can easily be carried out by hand for algebras of moderate $k$-dimensions,
say of dimensions $\le 50$ (see \cite{ZHpre}).

Due to the fact that the appearance of a principal left ideal $\Lambda p$ in
the second syzygy of a module $M$ can be related to the structure of
the first syzygy \cite{ZHinvent}, the above theorem, moreover, allows us
to build monomial relation algebras with prescribed finitistic dimensions
differing by at most 1.

Further examples resulting from this approach are as follows: There exist
monomial relation algebras for which the little finitistic dimension is not
attained on a cyclic module; more strongly, for each positive integer $b$,
there exists a monomial relation algebra whose little finitistic dimension is
not attained on a module of length at most $b$. On top of this, the mentioned
phenomena may depend on the choice of base field. See \cite{ZHfielddep}.
\medskip

{\bf Example for Method 3.} Let $\Gamma$ be the quiver
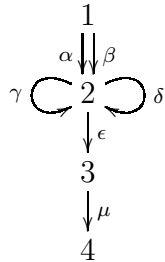
\ignore{
$$\xyoption{matrix}\xyoption{arrow}\xyoption{curve}\xymatrixrowsep{1.25pc}
\xy\xymatrix{
1 \arrow[d]<-0.4ex>_\alpha \arrow[d]<0.4ex>^\beta\\
2 \arrow[d]^\epsilon \lloopd{\gamma} \rloopd{\delta}\\
3 \arrow[d]^\mu\\
4}\endxy$$}
\noindent and let $\Lambda= k\Gamma/I$ be the monomial relation algebra whose
indecomposable projectives have the following graphs:

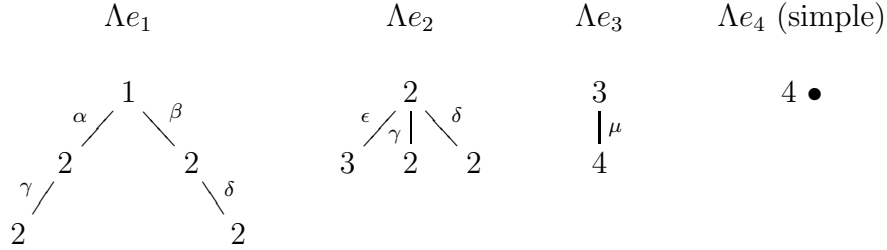
\vfill\eject
\vskip-0.25truein
\ignore{
\centerline{\hbox{
$$\xyoption{matrix}\xyoption{arrow}\xymatrixrowsep{1pc}\xymatrixcolsep{0.4pc}
\xy\xymatrix{
 &&\Lambda e_1 &&& &&&\Lambda e_2 &&&&\Lambda e_3 &&&\Lambda e_4\ {\rm
(simple)}\\
 &&1 \arrow@{-}[dl]_\alpha \arrow@{-}[dr]^\beta &&& &&&2
\arrow@{-}[dl]_\epsilon \arrow@{-}[d]_(0.6)\gamma \arrow@{-}[dr]^\delta &&&&3
\arrow@{-}[d]^\mu &&&4\ \bullet\\
 &2 \arrow@{-}[dl]_\gamma &&2 \arrow@{-}[dr]^\delta &&&&3 &2 &2 &&&4\\
2 &&&&2}\endxy$$}}}
\medskip

\noindent (Again we refer to \cite{ZHpre} for an interpretation of these
graphs.) In other words, the ideal
$I$ is generated by all those paths which do not make an appearance in any of
the above graphs. Using the method of
\cite{ZHpre}, one readily checks that the only paths $q$ with $\pdim \Lambda
q <\infty$ are $q=\epsilon$ and $q=\mu$, and that $\pdim \Lambda\epsilon =1$,
$\pdim\Lambda \mu =0$. Thus we obtain $s=1$ and conclude that
$$2\le \fd\Lambda \le \Fd\Lambda \le 3.$$
In fact, the little finitistic dimension of $\Lambda$ is 3, since $\pdim
\Lambda e_1/ \Lambda (\alpha+\beta) =3$.

\section{Method 4: Contravariant finiteness}

The following concept is due to Auslander and Smal\o\ \cite{AS} and was
further developed by Auslander-Buchweitz \cite{ABuch} in the commutative
case, by Auslander-Reiten \cite{AR91} in the context of finite dimensional
algebras. A full subcategory $\A$ of $\lamod$ is called {\it contravariantly
finite} if each module $M$ in $\lamod$ has an {\it $\A$-approximation} as
follows: there exists a homomorphism $f : A\rightarrow M$ with $A\in\A$ such
that each $g\in \Hom_\la (B,M)$ with $B\in\A$ factors through $f$, i.e.,
\ignore{
$$\xyoption{matrix}\xyoption{arrow}
\xy\xymatrix{
A \arrow[rr]^f &&M\\
 &&B \arrow[u]_g \arrow@{{}.>}[ull]^\exists}\endxy$$}

In the literature, the maps which we call $\A$-approximations are labeled
{\it right} $\A$-approximations, to account for the dual concept. It is well
known that, provided $M$ has an $\A$-approximation, the $\A$-approximations of
minimal length are isomorphic; it is thus unambiguous to speak of ``the'' {\it
minimal $\A$-approximation} of $M$ in that case. Recall that
$\pinf$ denotes the full subcategory of
$\lamod$ the objects of which are the modules of finite projective dimension.
As we mentioned in the synopsis of results above, contravariant finiteness of
$\pinf$ in $\lamod$ implies finiteness of $\fd \la$ \cite{AR91}. The major
advantage of this particular approach to classes of algebras satisfying the
second Finitistic Dimension Conjecture lies in the fact that, potentially, it
yields a great deal of additional information as a byproduct.  
\medskip

$\bullet$ If $\pinf$ is contravariantly finite in $\lamod$, and
$A_1,\dots,A_n$ are the minimal $\pinf$-approximations of the simple left
$\la$-modules, then $\fd\la =\max \{\pdim A_i \mid 1\le i\le n\}$. 
Moreover, an object in $\lamod$ has finite projective dimension if and only
if it is a direct summand of a module $X$ with a filtration
$X=X_0\supseteq X_1\supseteq \dots\supseteq X_d=0$ with $X_i/X_{i+1}\in
\{A_1,\dots,A_n\}$, up to isomorphism \cite{AR91}. 

In other words, the minimal
$\pinf$-approximations of the simples are the basic building blocks of the
modules of finite projective dimension in case $\pinf$ is contravariantly
finite. The drawback of this approach lies in the following catch: For few
classes of algebras is it known whether $\pinf$ is contravariantly finite
and, in general, the question is not at all easy to decide. Furthermore, even
when
$\pinf$ is known to be contravariantly finite, describing or computing the
minimal approximations of the simple modules is rather tricky. To date, the
following results have been obtained along this line:
\medskip

$\bullet$ If $\la$ is stably equivalent to a hereditary algebra, then $\pinf$
is contravariantly finite \cite{AR91}.

$\bullet$ Igusa-Smal\o-Todorov demonstrated with an example that $\pinf$ need
not be contravariantly finite, not even for a monomial relation algebra with
$J^3=0$ \cite{IST}. Their example will be presented below.

$\bullet$ If $\la$ is left serial, meaning that the indecomposable projective
left $\la$-modules are uniserial, then $\pinf$ is contravariantly finite, and
the minimal approximations of the simple left $\la$-modules are completely
understood \cite{BuHu}. In fact, given quiver and relations for $\la$, they
can be explicitly constructed. It turns out that, contrary to the widespread
belief that the left module theory over left serial algebras is fairly
simple, these minimal approximations attest to a rather sophisticated
structure in general. We illustrate the shape of the resulting approximations
with an example.
\medskip

{\bf Example for Method 4.} Let $\Gamma$ be the quiver
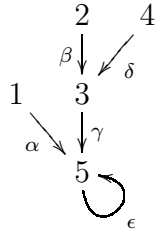
\ignore{
$$\xyoption{matrix}\xyoption{arrow}\xyoption{curve}
\xymatrixrowsep{1.25pc}\xymatrixcolsep{1pc}
\xy\xymatrix{
 &2 \arrow[d]_\beta &4 \arrow[dl]^\delta\\
1 \arrow[dr]_\alpha &3 \arrow[d]^\gamma\\
 &5 \ar@(d,r)_\epsilon
}\endxy$$}
\noindent and let $I\subseteq k\Gamma$ be such that the graphs of the
indecomposable projective left $\la$-modules are as follows:
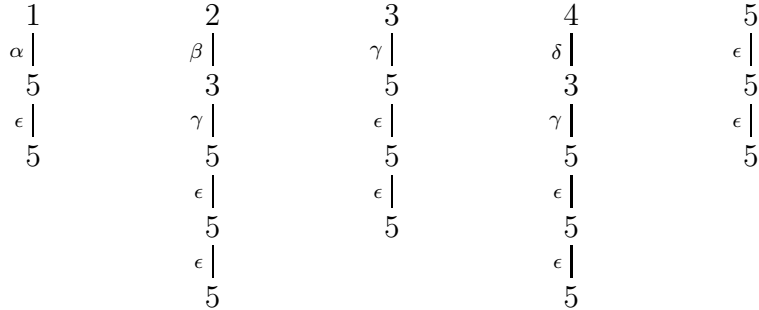
\ignore{
$$\xyoption{matrix}\xyoption{arrow}\xymatrixrowsep{1pc}
\xy\xymatrix{
1 \arrow@{-}[d]_\alpha &&2 \arrow@{-}[d]_\beta &&3 \arrow@{-}[d]_\gamma &&4
\arrow@{-}[d]_\delta &&5 \arrow@{-}[d]_\epsilon\\
5 \arrow@{-}[d]_\epsilon &&3 \arrow@{-}[d]_\gamma &&5 \arrow@{-}[d]_\epsilon
&&3 \arrow@{-}[d]_\gamma &&5 \arrow@{-}[d]_\epsilon\\
5 &&5 \arrow@{-}[d]_\epsilon &&5 \arrow@{-}[d]_\epsilon &&5
\arrow@{-}[d]_\epsilon &&5\\
 &&5 \arrow@{-}[d]_\epsilon &&5 &&5 \arrow@{-}[d]_\epsilon\\
 &&5 && &&5}\endxy$$}
\noindent Then the minimal $\pinf$-approximation $A_1$ of the simple $S_1$
going with the vertex 1 has graph
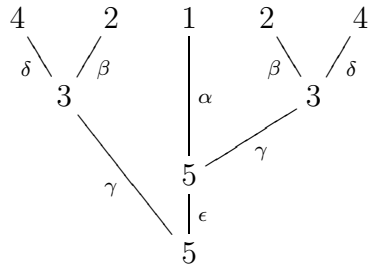
\ignore{
$$\xyoption{matrix}\xyoption{arrow}\xymatrixrowsep{1.25pc}\xymatrixcolsep{0.4pc}
\xy\xymatrix{
4 \arrow@{-}[dr]_\delta &&2 \arrow@{-}[dl]^\beta &&1 \arrow@{-}[dd]^\alpha
&&2 \arrow@{-}[dr]_\beta &&4 \arrow@{-}[dl]^\delta\\
 &3 \arrow@{-}[ddrrr]_\gamma &&& &&&3 \arrow@{-}[dlll]^\gamma\\
 &&&&5 \arrow@{-}[d]^\epsilon\\
 &&&&5}\endxy$$}
\medskip

$\bullet$ Happel and the author developed criteria for contravariant
finiteness and infinite dimensional substitutes for minimal
$\pinf$-approximations in case $\pinf$ fails to be contravariantly finite. To
define these substitutes, consider the following slight variant of
approximations as defined by Auslander-Smal\o: Suppose that $\C$ is a
subcategory of $\A$. A {\it $\C$-approximation of $M$ inside $\A$} is a map
$f \in \Hom_\la (A,M)$ with $A\in\A$ such that
\ignore{
$$\xyoption{matrix}\xyoption{arrow}
\xy\xymatrix{
A \arrow[rr]^f &&M\\
 &&C \arrow[u]_g \arrow@{{}.>}[ull]^\exists}\endxy$$}
\noindent for all $C\in\C$ and $g\in \Hom_\la(C,M)$. Note that, whenever $\C$
has finite representation type, $\C$-approximations of any module $M$ inside
$\A$ exist. The following ``phantoms'' are to take over the role of minimal
$\pinf$-approximations when the latter fail to exist. 

Fix $M\in\lamod$. A $\la$-module $H$ is called a {\it $\pinf$-phantom} of $M$
if there exists a nonempty finite subclass $\A$ of $\pinf$ such that each
$\mbox{add}(\A)$-approximation of $M$ inside $\pinf$ has $H$ as a subfactor;
direct limits of such modules $H$ are again labeled phantoms. In particular,
phantoms of the latter type need not be finitely generated. The crucial point
is that non-trivial phantoms always exist. Of course, the ``best'' phantom is
 the minimal $\pinf$-approximation of $M$ in the
case of existence. The remaining case is covered by the following theorem:
\medskip

$\bullet$ $M$ does not have a $\pinf$-approximation if and only if $M$ has
$\pinf$-phantoms of countably infinite $k$-dimension \cite{HaHu}. 
\medskip

This equivalence comes with instructions for finding phantoms over various
classes of algebras. We illustrate these concepts with several examples. The
first is a slight variant of the above-mentioned Igusa-Smal\o-Todorov
example. 
\medskip

{\bf Further Examples for Method 4.} (a) Let $\Gamma$ be the quiver
\ignore{
$$\xyoption{matrix}\xyoption{arrow}\xyoption{curve}\xymatrixcolsep{3pc}
\xy\xymatrix{
1 \arrow[r]<0.5ex>^\alpha \arrow[r]<-0.5ex>_\beta &2 \rloopd{\gamma}
}\endxy$$}
\noindent and suppose that the indecomposable projectives have the following
graphs:
\ignore{
$$\xyoption{matrix}\xyoption{arrow}\xymatrixrowsep{1pc}\xymatrixcolsep{0.4pc}
\xy\xymatrix{
 &&\Lambda e_1 &&&& &&&\Lambda e_2\\
 &&1 \arrow@{-}[dl]_\alpha \arrow@{-}[dr]^\beta &&&& &&&2
\arrow@{-}[d]^\gamma\\
 &2 \arrow@{-}[dl]_\gamma &&2 &&& &&&2\\
2}\endxy$$}
\noindent Then the simple left $\la$-module $S_1$ fails to have a
$\pinf$-approximation. The following graph represents an infinite dimensional
$\pinf$-phantom $H_1$ of $S_1$:
\ignore{
$$\xyoption{matrix}\xyoption{arrow}\xymatrixrowsep{1.25pc}\xymatrixcolsep{0.4pc}
\xy\xymatrix{
1 \arrow@{-}[dr]_\beta &&1 \arrow@{-}[dr]_\beta &&1 \arrow@{-}[dr]_\beta &&
\arrow@{{}.{}}[rr] &&\\
 &2 \arrow@{-}[ur]^\alpha &&2 \arrow@{-}[ur]^\alpha &&2 & \arrow@{{}.{}}[rr]
&&}\endxy$$}
\noindent This example might leave the impression that the existence of a
module $H_1$ of finite projective dimension, with graph as above, already
prevents
$S_1$ from having a $\pinf$-approximation. The next example counters this
impression.

(b) This time, the quiver $\Gamma$ is 
\ignore{
$$\xyoption{matrix}\xyoption{arrow}\xyoption{curve}\xymatrixcolsep{3pc}
\xy\xymatrix{
3 \arrow[r]^\delta &1 \arrow[r]<0.5ex>^\alpha \arrow[r]<-0.5ex>_\beta &2
\rloopd{\gamma} }\endxy$$}
\noindent and the indecomposable projectives in $\lamod$ have the graphs 
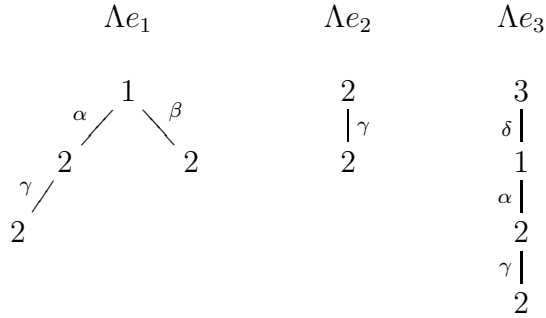
\ignore{
$$\xyoption{matrix}\xyoption{arrow}\xymatrixrowsep{1pc}\xymatrixcolsep{0.4pc}
\xy\xymatrix{
 &&\Lambda e_1 &&&&&\Lambda e_2 &&&&\Lambda e_3\\
 &&1 \arrow@{-}[dl]_\alpha \arrow@{-}[dr]^\beta &&&&&2
\arrow@{-}[d]^\gamma &&&&3 \arrow@{-}[d]_\delta\\
 &2 \arrow@{-}[dl]_\gamma &&2 &&&&2 &&&&1 \arrow@{-}[d]_\alpha\\
2 &&&& &&&& &&&2 \arrow@{-}[d]_\gamma\\
 &&&& &&&& &&&2}\endxy$$}
\noindent Note that the graphs of $\la e_1$ and $\la e_2$ are as under (a).
In particular, there is again a unique infinite dimensional module $H_1$ of
finite projective dimension having the graph displayed above. However, this
time,
$H_1$ is not a
$\pinf$-phantom of
$S_1$. In fact, $\pinf$ is contravariantly finite in this example, the
minimal $\pinf$-approximation of $S_1$ being of the form
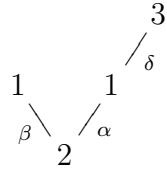
\ignore{
$$\xyoption{matrix}\xyoption{arrow}\xymatrixrowsep{1pc}\xymatrixcolsep{0.4pc}
\xy\xymatrix{
 &&&3 \arrow@{-}[dl]^\delta\\
1 \arrow@{-}[dr]_\beta &&1 \arrow@{-}[dl]^\alpha\\
 &2}\endxy$$}
\noindent Our final example is to indicate how instable both of the
conditions, contravariant finiteness of $\pinf$ and failure thereof, are in
general. 

(c) Take over the quiver of example (b), as well as the graphs of the
indecomposable projectives $\la e_1$ and $\la e_2$. Only the graph of $\la
e_3$ is modified, through the deletion of one relation. 
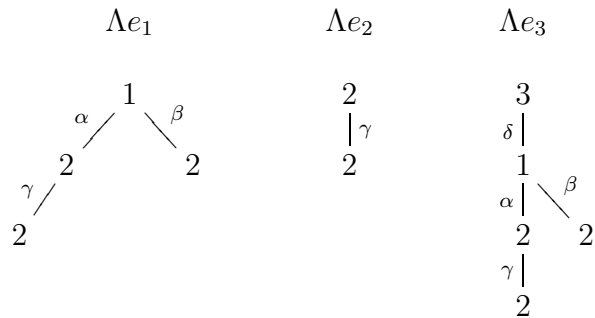
\ignore{
$$\xyoption{matrix}\xyoption{arrow}\xymatrixrowsep{1pc}\xymatrixcolsep{0.4pc}
\xy\xymatrix{
 &&\Lambda e_1 &&&&&\Lambda e_2 &&&&\Lambda e_3\\
 &&1 \arrow@{-}[dl]_\alpha \arrow@{-}[dr]^\beta &&&&&2
\arrow@{-}[d]^\gamma &&&&3 \arrow@{-}[d]_\delta\\
 &2 \arrow@{-}[dl]_\gamma &&2 &&&&2 &&&&1 \arrow@{-}[d]_\alpha
\arrow@{-}[dr]^\beta\\ 
2 &&&& &&&& &&&2 \arrow@{-}[d]_\gamma &2\\
 &&&& &&&& &&&2}\endxy$$}
\noindent We display two $\pinf$-phantoms of $S_1$ having infinite
$k$-dimension:
\ignore{
$$\xyoption{matrix}\xyoption{arrow}\xymatrixrowsep{1.25pc}\xymatrixcolsep{0.4pc}
\xy\xymatrix{
 &&&& &&&& &&&& &&&&3 \arrow@{-}[d]_\delta &&3 \arrow@{-}[d]_\delta &&
\arrow@{{}.{}}[rr] &&\\
1 \arrow@{-}[dr]_\beta &&1 \arrow@{-}[dr]_\beta &&1 \arrow@{-}[dr]_\beta &&
\arrow@{{}.{}}[rr] &&  &&&{\rm and} &&&1 \arrow@{-}[dr]_\beta &&1
\arrow@{-}[dr]_\beta &&1 \arrow@{-}[dr]_\beta &&
\arrow@{{}.{}}[rr] &&\\
 &2 \arrow@{-}[ur]^\alpha &&2 \arrow@{-}[ur]^\alpha &&2 & \arrow@{{}.{}}[rr]
&& &&& &&& &2 \arrow@{-}[ur]^\alpha &&2 \arrow@{-}[ur]^\alpha &&2 &
\arrow@{{}.{}}[rr] &&}\endxy$$}
\noindent In particular, these phantoms show that $S_1$ does not have a
$\pinf$-approximation in this example.

\end{document}